\documentstyle[a4,12pt,amssymb]{article}
\newtheorem{theorem}{Theorem}
\newtheorem{lemma}{Lemma}
\newtheorem{prop}{Proposition}
\newcommand{\ra}{\longrightarrow}
\newcommand{\da}{\downarrow}
\newcommand{\cO}{\hat {\cal O}}
\newcommand{\cF}{\hat {\cal F}}
\newcommand{\pr}{\prod\limits}
\def\R{{\mathbb R}}
\def\Z{{\mathbb Z}}

\def\C{{\mathbb C}}

\def\A{{\mathbb A}}
\author{}
\date{}
\begin{document}

\vspace*{2cm}
\begin{center}
{\bf KRICHEVER CORRESPONDENCE \\
FOR ALGEBRAIC SURFACES} \\
\par\bigskip
{\sc A. N. PARSHIN} \\
\par\medskip
Steklov Mathematical Insttiute, Moscow \\
Institut f\"ur Algebra und Zahlentheorie der \\
Technische Universit\"at Braunschweig \\
\par\bigskip
Bericht 99/17
\end{center}
\par\bigskip
In 70's there was discovered a construction how to attach to some
algebraic-geometric data an infinite-dimensional subspace in the space
$k((z))$ of the Laurent power series. The construction was
successfully used in the theory of integrable systems, particularly, for
the KP and KdV equations \cite{K, SW}.
There were also found some applications to the moduli of algebraic curves
\cite{ADKP, BS}.
Now it is known as the Krichever correspondence or the Krichever map
\cite{ADKP, M, AMP, Q, BF}. The original work by I. M. Krichever has also
included commutative rings of differential operators as a third part of
the correspondence.

The map we want to
study here was first described in an explicit way by G. Segal and G. Wilson
\cite{SW}. They have used an analytical version of the infinite dimensional
Grassmanian introduced by M. Sato \cite{SS, PS}. In the sequel we consider a
purely algebraic approach as developed in \cite{M}.

Let us just note that
the core of the construction is an embedding of the affine coordinate ring
on an algebraic  curve into the field $k((z))$ corresponding to the power
decompositions in a point at infinity
(the details see below in section 2). In number theory this corresponds to
an embedding of the ring of algebraic integers to the fields $\C$ or $\R$.
The latter one is well known starting from the XIX-th century. The  idea
introduced by Krichever was to insert the local parameter $z$. This
trick looking so simple enormously extends the area of the correspondence.
It allows to consider all algebraic curves simultaneously.

But there still remained a hard restriction by the case of curves, so by
dimension 1. Recently, it was pointed out by the author \cite{P3} that there
are some connections between the theory of the KP-equations and the theory of
$n$-dimensional local fields \cite{P1}, \cite{FP}. From this point of view
it becomes clear that the Krichever construction should have a generalization
to the case of higher dimensions. This generalization is suggested in the
paper  for the case of algebraic surfaces (see theorem 3 in section 2). A
further generalization  to the case of arbitrary dimension was recenly
proposed by D. V. Osipov \cite{O}.

Let us also note that the construction of the restricted adelic complex
in section 1 is of an independent interest, also in arithmetics. It has
already appeared in a description of vector bundles on algebraic surfaces
\cite{P2}.
\par\medskip
This work was partly done during my visit to the
Institut f\"ur Algebra
und Zahlentheorie der Technische Universit\"at Braunschweig in January 1999.
I am very much grateful to Professor Hans Opolka for the hospitality.

\section{Adelic Complexes}
We first discuss the adelic complexes for the case of dimension 1.
Concerning a definition of the adelic notions we refer to \cite{FP},\cite{H}.
We also note that the sign $\prod$ denotes the adelic product.

Let $C$ be an projective algebraic curve over a field $k$, $P$ be
a smooth point  and $\eta$ a general point on $C$. Let ${\cal F}$ be a
torsion free coherent sheaf on $C$.
\begin{prop}. The following complexes are quasi-isomorphic:
\begin{quotation}
i) adelic complex
$$ {\cal F}_{\eta} \oplus \pr_{x \in C} \cF_x \ra \prod_{x \in C}
(\cF_x \otimes_{\cO_x} K_x)
$$
ii) the complex
$$ W \oplus \cF_P \ra \cF_P \otimes_{\cO_P} K_P
$$
where  $W = \Gamma(C - P, {\cal F}) \subset \cF_{\eta}$.
\end{quotation}
\end{prop}
{\sc Proof} will be done in two steps. First, the adelic complex
contains a trivial exact subcomplex
$$  \pr_{x \in U} \cF_x \ra \pr_{x \in U} \cF_x,
$$
where $U = C - P$. The quotient-complex is equal to
$$
{\cal F}_{\eta} \oplus \cF_P \ra \prod_{x \in U}
(\cF_x \otimes K_x)/\cF_x  \oplus \cF_P.
$$
It has a surjective homomorphism to the exact complex
$$
{\cal F}_{\eta}/W \ra \prod_{x \in U}
(\cF_x \otimes K_x)/\cF_x.
$$
The exactness of the complex is the strong approximation theorem
for the curve $C$ \cite{B}[ch.II, \S 3, corollary of prop. 9; ch. VII,
\S, prop. 2]. The kernel of this
surjection will be the second complex from proposition.
\par\smallskip

Now we go to the case of dimension 2.
Let $X$ be a projective irreducible algebraic surface over a field
     $k$, $C \subset X$ be an irreducible projective curve, and
     $P \in C$ be a smooth point on both $C$ and $X$. Let ${\cal F}$
be a torsioh free coherent sheaf on $X$.


{\sc Definition 1}. Let $x \in C$. We let
     $$ B_x({\cal F}) = \bigcap\limits_{D \ne C}((\cF_x \otimes K_x) \cap
(\cF_x \otimes \cO_{x,D})),$$
     where the intersection is done inside the group $\cF_x \otimes K_{x}$,
     $$ B_C({\cal F}) = (\cF_C \otimes K_C) \cap
     (\bigcap\limits_{x \ne P} B_x),$$
     where the intersection is done inside $\cF_x \otimes K_{x,C}$,
     $$ A_C({\cal F}) = B_C({\cal F}) \cap \cF_C,$$
     $$ A({\cal F}) = \cF_{\eta} \cap (\bigcap\limits_{x \in X-C} \cF_x).$$

We will freely use the following shortcuts:
$$
\begin{array}{lll}
K\cF_x & = &  \cF_x \otimes_{\cO_x} K_x, \\
K\cF_D & = &  \cF_D \otimes_{\cO_D} K_D, \\
\cF_{x,D} & = &  \cF_x \otimes_{\cO_x} {\cal O}_{x,D}, \\
K\cF_{x,D} & = &  \cF_x \otimes_{\cO_x} K_{x,D}.
\end{array}
$$
Next, we need two lemmas connecting the adelic complexes on $X$ and $C$.
They are the versions of the relative exact sequences, see \cite{P1},
\cite{FP}.
The curve $C$ defines the following ideals:
$$
K_{x,C} \supset \cO_{x,C} \dots \supset \wp^n_{x,C} \supset \dots,
$$
$$
K_C \supset \cO_C \dots \supset \wp^n_C \supset \wp^{n+1}_C \supset \dots,
$$
$$
K_x \supset \cO_x \dots \supset \wp^n_x \dots ,
$$
and $\wp_x = \cO_x \cap \wp_{x,C}$.

\begin{lemma}. We assume that the curve $C$ is a locally complete intersection.
Let $N_{X/C}$ be the normal sheaf for the curve $C$  in $X$.
For all $n \in \Z$ the maps
$$
\prod_{x \in C} \wp^n_{x,C}\cF_{x,C}/\wp^{n+1}_{x,C}\cF_{x,C}   \ra
\A_{C,01}({\cal F} \otimes \check{N}_{X/C}^{\otimes n}),
$$
$$
\prod_{x,C} \wp^n_x\cF_{x}/\wp^{n+1}_x\cF_{x} \ra \A_{C,1}({\cal F}
\otimes \check{N}_{X/C}^{\otimes n}),
$$
$$
\wp^n_C\cF_{C}/\wp^{n+1}_C\cF_{C} \ra \A_{C,0}({\cal F}
\otimes \check{N}_{X/C}^{\otimes n}),
$$
are bijective.
\end{lemma}
In general, we have an exact sequence
$$
0 \ra  {\cal J}^{n+1}  \ra {\cal J}^{n} \ra {\cal J}^{n}\vert_C \ra 0
$$
where ${\cal J} \subset  {\cal O}_X$ is an ideal defining the curve $C$.
In our case ${\cal J} = {\cal O}_X(-C)$ and $N_{X/C} = {\cal O}_X(C)\vert_C$.
Thus the isomorphisms from the lemma are coming from the
exact relative sequence
$$
0 \ra \A_X({\cal F}(-(n+1)C)) \ra \A_X({\cal F}(-nC)) \ra
\A_C({\cal F}(-nC)\vert_C) \ra 0.
$$

\begin{lemma}. Let $P \in C$. For all $n \in {\bf Z}$  the complex
$$
\wp^n_C\cF_{C}/\wp^{n+1}_C\cF_{C} \oplus \prod_{x \in C}
\wp^n_x\cF_{x}/\wp^{n+1}_x\cF_{x} \ra
\prod_{x \in C} \wp^n_{x,C}\cF_{x,C}/\wp^{n+1}_{x,C}\cF_{x,C}
$$
is quasi-isomorphic to the complex
$$
(A_C({\cal F}) \cap \wp^n_C\cF_{C})/(A_C({\cal F}) \cap \wp^{n+1}_C\cF_{C})
\oplus \wp^n_P\cF_{P}/\wp^{n+1}_P\cF_{P} \ra
\wp^n_{P,C}\cF_{P,C}/\wp^{n+1}_{P,C}\cF_{P,C}.
$$
\end{lemma}
This lemma is an extension of the proposition 1 above. The proves of the
both lemmas are straightforward and we will skip them.

\begin{theorem}.
Let $X$ be a projective irreducible algebraic surface over a field
     $k$, $C \subset X$ be an irreducible projective curve, and
     $P \in C$ be a smooth point on both $C$ and $X$. Let ${\cal F}$
be a torsioh free coherent sheaf on $X$.

Assume that the the surface $X - C$ is affine.
Then the following complexes are quasi-isomorphic:
\begin{quotation}

i)  the adelic complex
$$
\cF_{\eta}  \oplus  \pr_{D} \cF_D   \oplus  \pr_{x} \cF_x
 \ra
\prod_{D} (\cF_D \otimes K_D) \oplus  \prod_{x} (\cF_x \otimes K_x) \oplus
\prod_{x \in D} (\cF_x \otimes \cO_{x,D})
 \ra
$$
$$
\ra \prod_{x \in D} (\cF_x \otimes K_{x,D})
$$
for the sheaf ${\cal F}$ and

ii) the complex
$$
A({\cal F}) \oplus A_C({\cal F}) \oplus \cF_P \ra B_C({\cal F}) \oplus
B_P({\cal F}) \oplus (\cF_P \otimes
\cO_{P,C}) \ra \cF_P \otimes  K_{P,C}
$$

\end{quotation}
\end{theorem}

{\sc  Proof} will be divided into several steps. We will subsequently
     transform the adelic complex checking that every time we get a
     quasi-isomorphic complex.

     {\bf Step I}.  Consider the diagram
{\footnotesize
$$
\begin{array}{cccccccccc}
&& \pr_{D \ne C} \cF_D & \oplus & \pr_{x \in U} \cF_x &
\ra & \pr_{D \ne C} \cF_D  & \oplus & \pr_{x \in U} \cF_x & \oplus \\

&& \da && \da && \da && \da & \\

\cF_{\eta} & \oplus & \pr_{D} \cF_D  & \oplus & \pr_{x} \cF_x
& \ra &
\prod_{D} K\cF_{D} & \oplus & \prod_{x} K\cF_{x} & \oplus \\

\parallel && \da && \da && \da && \da & \\

\cF_{\eta} & \oplus &  \cF_C  & \oplus & \pr_{x \in C} \cF_x
& \ra &
(\prod_{D \ne C} K\cF_{D}/\cF_D \oplus K\cF_{C}) & \oplus &
(\prod_{x \in U} K\cF_{x}/\cF_x \oplus \prod_{x \in C} K\cF_{x}) &
\oplus
\end{array}
$$
$$
\begin{array}{cccccc}
\oplus & \prod_{x \in D \ne C} \cF_{x,D} &
\ra & \prod_{x \in D \ne C} \cF_{x,D} && \\
& \da & &  \da &&\\
\oplus & \prod_{x \in D} \cF_{x,D} &
\ra & \prod_{x \in D} K\cF_{x,D} && \\
 & \da && \da &&\\

\oplus & \prod_{x \in C} \cF_{x,C} &
\ra & \prod_{x \in D \ne C} K\cF_{x,D}/\cF_{x,D} & \oplus &
\prod_{x,C} K\cF_{x,C}
\end{array}
$$}
where $U = X - C$. The middle
     row is the full adelic complex and the first row is an exact
     subcomplex. The commutativity of the upper squares is obvious.
     The exactness follows from the  trivial
\begin{lemma}. Let $f_{1,2}: A_{1,2} \ra B$ be  homomorphisms of abelian
groups.
     The complex
     $$ 0 \ra A_1 \oplus A_2 \ra A_1 \oplus A_2 \oplus B \ra B \ra 0, $$
     where $ (a_1 \oplus a_2) \mapsto (a_1 \oplus -a_2 \oplus -f(a_1) +f(a_2)),
 (a_1 \oplus a_2 \oplus b) \mapsto (f(a_1) + f(a_2) + b)$,
     is exact.
\end{lemma}
The third row in the diagram is a quotient-complex by the subcomplex and
     we conclude that it is quasi-isomorphic to the adelic complex.

		 {\bf Step II}. We can make the same step with the adelic complex
for the sheaf ${\cal F}$ on the surface $U$. By assumption the surface
     $U$ is affine  and we get an {\em exact} complex
$$
\cF_{\eta}/A \ra \prod_{D \ne C} (\cF_D \otimes K_D)/\cF_D \oplus
\prod_{x \in U} (\cF_x \otimes K_x)/\cF_x
     \ra
$$
$$
  \prod_{\begin{array}{c}
               x \in U \\
              x \in D \ne C
               \end{array}}
     (\cF_x \otimes K_{x,D})/(\cF_x \otimes \cO_{x,D}),
$$
where $A = \Gamma(U, {\cal F})$.
\begin{lemma}. The complex
$$
0 \ra \prod_{x \in C} (\cF_x \otimes K_x)/B_x({\cal F})
     \ra \prod_{\begin{array}{c}
               x \in C \\
              x \in D \ne C
               \end{array}}
     (\cF_x \otimes K_{x,D})/(\cF_x \otimes \cO_{x,D}) \ra 0
$$
is exact.
\end{lemma}
{\sc Proof}. The injectivity follows directly from the definition of
     the ring $B_x$. The surjectivity is the local strong approximation
     around the point $X \in X$ (see \cite{P1}[\S 1],\cite{FP}[ch.4]).

     {\bf Step III}.  Take the sum of the  two complexes from
     step II. Then we have a map of the complex we got in the step I to
     this complex
{\footnotesize
$$
\begin{array}{cccccccccc}

\cF_{\eta} & \oplus &  \cF_C  & \oplus & \pr_{x} \cF_x
& \ra &
(\prod_{D \ne C} K\cF_{D}/\cF_D \oplus K\cF_{C}) & \oplus &
(\prod_{x \in U} K\cF_{x}/\cF_x  \oplus  \prod_{x \in C} K\cF_{x})
& \oplus \\

\da && \da && \da && \da && \da &  \\

\cF_{\eta}/A & \oplus & (0) & \oplus & (0)
& \ra &
\prod_{D \ne C} K\cF_{D}/\cF_D & \oplus &
(\prod_{x \in U} K\cF_{x}/\cF_x  \oplus  \prod_{x \in C} K\cF_{x}/B_x)
& \oplus
\end{array}
$$
$$
\begin{array}{cccccc}

\oplus & \prod_{x \in C} \cF_{x,C} &
\ra & \prod_{x \in D \ne C} K\cF_{x,D}/\cF_{x,D}  &
\oplus & \prod_{x \in C} K\cF_{x,C}  \\

& \da && \da && \da  \\

\oplus & (0) &
\ra & \prod_{\begin{array}{c}
x \in U \\
x \in D \\
D \ne C
\end{array}} K\cF_{x,D}/\cF_{x,D}
 \oplus
\prod_{\begin{array}{c}
x \in C \\
x \in D \\
D \ne C
\end{array}} K\cF_{x,D}/\cF_{x,D} & \oplus & (0)
\end{array}
$$}
For this map all the
     components which do not have arrows are mapped to zero. The diagram
     is commutative and the kernel of the map is equal to
$$
A \oplus \cF_C \oplus \pr_{x \in C} \cF_x  \ra
     K\cF_C \oplus \prod_{x \in C} B_x({\cal F}) \oplus
\prod_{x \in C} K\cF_x  \ra \prod_{x \in C} K\cF_{x,C}.
$$
We conclude that this complex is quasi-isomorphic to the adelic complex.

     {\bf Step IV}. Using the embedding $\cF_x \ra B_x({\cal F})$ and lemma 3
     we have an exact complex and it's embedding into the complex of the
     step III:
{\small
$$
\begin{array}{cccccccccccc}
&&&& \pr_{x \in C-P} \cF_x & \ra &&& \pr_{x \in C-P} B_x({\cal F}) & \oplus &
\pr_{x \in C-P} \cF_x &  \ra \\

&&&& \da &&&& \da && \da &  \\

A & \oplus & \cF_C & \oplus & \pr_{x \in C}  \cF_x & \ra &
K\cF_C  & \oplus & \prod_{x \in C} B_x({\cal F}) & \oplus &
\prod_{x \in C} \cF_{x,C}
  &   \ra
\end{array}
$$
$$
\begin{array}{cc}
\ra & \pr_{x \in C-P} B_x({\cal F})  \\
& \da \\
\ra &  \prod_{x \in C} K\cF_{x,C}
\end{array}
$$
}

As a result we get the factor-complex
$$
A \oplus \cF_C \oplus \cF_P \ra K\cF_C \oplus B_P({\cal F}) \oplus
     \prod_{x \in C-P} \cF_{x,C}/\cF_x \oplus \cF_{P,C} \ra
$$
$$
\ra  \prod_{x \in C-P} K\cF_{x,C}/B_x({\cal F}) \oplus K\cF_{P,C}.
$$

     {\bf Step V}. Now we need
\begin{lemma} The complex
$$
0 \ra  (\cF_C \otimes K_C)/B_C({\cal F})
     \ra \prod_{x \in C-P}
     (\cF_x \otimes K_{x,C})/B_x({\cal F})  \ra 0
$$
is exact.
\end{lemma}
{\sc Proof}. The injectivity is again the definition of the $B_C$ and
     the surjectivity follows from the strong approximation on the
     curve $C$ (\cite{B}) and  lemma 2 above.

     As a corollary we have an isomorphism
$$
\cF_C/A_C({\cal F}) \stackrel{\cong}{\ra} \prod_{x \in C-P} (\cF_x \otimes
\cO_{x,C})/\cF_x,
$$
where
$$ A_C({\cal F}) := B_C({\cal F}) \cap \cF_C. $$

Combining the isomorphisms from the lemma and its corollary into a
     single complex of length 2, we get the diagram
{\small
$$
\begin{array}{cccccccccccc}
A & \oplus & \cF_C & \oplus & \cF_P
& \ra &
K\cF_{C} & \oplus & B_P({\cal F}) & \oplus &
     (\prod_{x \in C-P} \cF_{x,C}/\cF_x \oplus \cF_{P,C})
& \ra \\

 \da && \da && \da && \da && \da && \da & \\

(0) & \oplus &   \cF_C/A_C  & \oplus & (0)
& \ra &
K\cF_C/B_C & \oplus &  (0) & \oplus & \prod_{x \in C-P} \cF_{x,C}/\cF_x
     & \ra
\end{array}
$$
$$
\begin{array}{cccc}
\ra  &
\prod_{x \in C-P} K\cF_{x,C}/B_x({\cal F}) & \oplus & K\cF_{P,C} \\
 & \da && \da \\
\ra  & \prod_{x \in C-P} K\cF_{x,C}/B_x({\cal F}) & \oplus & (0)
\end{array}
$$
}
The kernel of the map of the complexes is obviously equal to

$$
A({\cal F}) \oplus A_C({\cal F}) \oplus \cF_P \ra B_C({\cal F}) \oplus
B_P({\cal F}) \oplus (\cF_P \otimes
\cO_{P,C}) \ra  (\cF_P \otimes K_{P,C})
$$
and we arrive to the conclusion of the theorem.

{\sc Remark 1}. Sometimes we will call the complex from the theorem as the
{\em restricted} adelic complex.

\begin{lemma}.Let $X$ be a projective irreducible variety over a field $k$
and ${\cal O}(1)$ be a very ample sheaf on $X$.
Then
\begin{enumerate}
\item The following conditions are equivalent
     \begin{quotation}
    i) $X$ is a Cohen-Macaulay variety

    ii) for any locally free sheaf ${\cal F}$ on $X$ and $i < \mbox{dim}(X)$
     ~$H^{i}(X, {\cal F}(n)) = (0)$ for $n << 0$
     \end{quotation}
\item If $X$ is normal of dimension $> 1$ then for any locally free sheaf
${\cal F}$ on $X$~$H^{1}(X, {\cal F}(n)) = (0)$ for $n << 0$
\end{enumerate}
\end{lemma}
{\sc Proof} see in \cite{Har}[ch. III, Thm. 7.6, Cor. 7.8]. We only note that
the last statement is known as the lemma of Enriques-Severi-Zariski. For
dimension 2 every normal variety is Cohen-Macaulay and thus the second claim
follows from the first one.

\begin{prop}. Let ${\cal F}$ be a locally free  coherent sheaf on the
projective irreducible surface $X$.

Assume that the local rings of the $X$ are Cohen-Macaulay and the curve
$C$ is a locally complete intersection. Then, inside the
field $K_{P,C}$, we have
$$B_C({\cal F}) \cap B_P({\cal F}) = A({\cal F}).$$
\end{prop}
{\sc Proof} will be done in several steps.

{\sc Step 1}. If we know the proposition for a sheaf ${\cal F}$ then it
is true for the sheaf ${\cal F}(nC)$ for any $n \in \Z$. Thus taking
a twist by  ${\cal O}(n)$ we can assume that $\mbox{deg}_C({\cal F}) < 0$.

{\sc Step 2}.
Now we show that $A_C({\cal F}) \cap \cF_P = (0)$. The filtrations from
lemma 1 gives the corresponding filtration of the group $A_C({\cal F})$.
Lemma 2 implies that
$$ \frac{(A_C({\cal F}) \cap \cF_P) \cap \wp^n\cF_P}{(A_C({\cal F}) \cap
\cF_P) \cap \wp^{n+1}\cF_P}
\cong \Gamma(C, {\cal F} \otimes \check{N}^{\otimes n}_{X/C}).
$$
Since $\mbox{deg}_C({\cal F}) < 0$,~$N_{X/C} = {\cal O}_X(C)\vert_C$ and
$\mbox{deg}_C(N_{X/C}) > 0$
we get that the last group is trivial.

{\sc Step 3}.
The next step is to prove the equality:
$$
B_C({\cal F}(-D)) \cap B_P({\cal F}(-D)) = A({\cal F}(-D)),
$$
where $D$ is an sufficiently ample divisor on $X$ distinct from the curve
$C$. By theorem 1 the cohomology of ${\cal F}_X(-D)$
can be computed from the complex
$$
A({\cal F}(-D)) \oplus A_C({\cal F}(-D)) \oplus \cF_P(-D) \ra
B_C({\cal F}(-D)) \oplus B_P({\cal F}(-D)) \oplus
\cF_{P,C}(-D)
$$
$$
\ra   K{\cal F}_{P,C}.
$$
Now take $a_{01} \in B_C({\cal F}(-D)), a_{02} \in B_P({\cal F}(-D))$ such
that $a_{01} + a_{02} = 0$.
They define an element $(a_{01} \oplus a_{02} \oplus 0)$ in the middle
component of the complex.
By our condition for $D$ and the lemma 6 we have $H^1(X, {\cal F}_X(-D)) = (0)$
and thus there exist $a_0 \in A({\cal F}(-D)), a_1 \in A_C{\cal F}(-D),
a_2 \in \cF_P(-D)$ such that $a_{01} = a_0 - a_1, a_{02} = a_2 - a_0, 0 =
a_1 - a_2$.

By the second step $a_1 = a_2 \in (A_C({\cal F}(-D)) \cap \cF_P(-D)) \subset
A_C({\cal F}) \cap \cF_P = (0)$ and, consequently, we have
$a_{01} (=-a_{02}) \in A({\cal F}(-D))$.

{\sc Step 4}.
The last step is to take two distinct divisors $D, D'$ such that
$D \cap D' \subset C$. Since $C$ is a hyperplane section  we can choose for
$D, D'$ two hyperplane sections whose intersection belongs to $C$. Therefore
their ideals in the ring $A({\cal F})$ are relatively prime and
$$ A({\cal F}) = A({\cal F}(-D)) + A({\cal F}(-D')) \ni 1 = a + a' , a
\in A({\cal F}(-D)), a' \in A({\cal F}(-D')).$$
If now $b \in B_C({\cal F}) \cap B_P({\cal F})$, then $b = ba + ba'$,
where $ba \in B_C({\cal F}(-D)) \cap B_P({\cal F}(-D)), ba' \in
B_C({\cal F}(-D')) \cap B_P({\cal F}(-D'))$. We see that $b \in A({\cal F})$
by the previous step.

{\sc Remark 2}.The method we have used  cannot be applied if our
variety is not Cohen-Macaulay (by lemma 6 above). It would be interesting to
know how to extend
the result to the arbitrary surfaces $X$ and the sheaves ${\cal F}$ such
that ${\cal F}$ are locally free outside $C$. The last condition is really
necessary.


{\sc Remark 3}. This proposition is a version for the reduced adelic complex
of the corresponding result for the full complex. Namely, $\A_{X,01} \cap
\A_{X,02} = \A_{X,0}$, see \cite[ch.IV]{FP}. This should be generalized
to arbitrary dimension $n$ in the following way.

Let $I, J \subset [0, 1, \dots, n]$ and
$$
\A_{X, I}({\cal F}) = (\prod_{\{codim \eta_0, codim \eta_1, \dots \} \in I}
K_{\eta_0,\eta_1,\dots}) \otimes {\cal F}_{\eta_0}) ~\bigcap ~\A_X({\cal F}).
$$
Then we have
$$
\A_{X, I}({\cal F}) \cap \A_{X, J}({\cal F}) = \A_{X, I \cap J}({\cal F})
$$
for a locally free ${\cal F}$ and a Cohen-Macaulay $X$.
\par\medskip
{\sc Example}. Let $X = {\bf P}_2 \supset C = {\bf P}_1 \ni P$. We introduce
homogenous coordinates $(x_0:x_1:x_2)$ such that  $C = (x_0 = 0); P =
(x_0 = x_1 = 0$ and
$U = X - C = Spec k[x,y]$ with $x = x_1/x_0, y = x_2/x_0$.
Then $k(C) = k(y/x), x^{-1}$ is the last parameter for any two-dimensional
local field $K_{Q,C}$ with $Q \ne P$. For local field $K_{P,C}$ we have
$$ K_{P,C} = k((u))((t)), u = xy^{-1}, t = y^{-1}.$$

Then we can easily compute all the rings from the complex of theorem 1 for
the sheaf ${\cal O}_X$.
$$
\begin{array}{ccccl}
&& B_P & = & k[[u]]((t)) \\
&& B_C & = & k[u^{-1}]((u^{-1}t)) \\
&& \cO_{P,C} & = & k((u))[[t]] \\
A & = & \Gamma(U, {\cal O}_X) & = & k[ut^{-1}, t^{-1}] \\
&& A_C & = & k[u^{-1}][[u^{-1}t]] \\
&& \cO_P & = & k[[u,t]]
\end{array}
$$
We can draw the subspaces as some subsets of the plane according to the
supports of the elements of the subspaces (on the plane with
coordinates $(i,j)$ for elements $u^it^j \in K_{P,C}$. Then the first
three subspaces $B_P, B_C, \cO_{P,C}$ will correspond to some halfplanes
and the subspaces $A, A_C, \cO_P$ to the intersections of them.

\section{Main Theorem}

We need the following well known result.
\begin{lemma}. Let $X$ be an projective variety, ${\cal F}$ be a coherent
sheaf on õXõ and $C$ be an ample divisor on $X$. If
$$
S = \oplus_{n \ge 0} \Gamma(X, {\cal O}_X(nC)),~F = \oplus_{n \ge 0}
\Gamma(X, {\cal F}(nC)),
$$
then
$$ X \cong Proj(S),~{\cal F} \cong Proj(F).$$
\end{lemma}
{\sc Proof}. Let $mC$ be a very ample divisor,$S = \oplus_{n \ge 0} S_n$ and
$S': = \oplus_{n \ge 0} S_{nm}$. Then by \cite{GD}[prop. 2.4.7]
$$Proj(S') \cong Proj(S).$$

The divisor $mC$ defines an embedding $i: X \ra {\bf P}$ to a projective
space such that $i^*{\cal O}_{\bf P}(1) = {\cal O}_X(mC)$. Let ${\cal J}_X \subset
{\cal O}_{\bf P}$ be an ideal defined by $X$. If

$$ I: = \oplus_{n \ge 0} \Gamma({\bf P}, {\cal J}_X(n)), $$
$$ A: = \oplus_{n \ge 0} \Gamma({\bf P}, {\cal O}_{\bf P}(n)), $$
then $I \subset A$ and by \cite{GD}[prop. 2.9.2]
	$$ Proj(A/I) \cong X.$$

	We have an exact sequence of sheaves
	$$0 \ra {\cal J}_X(n) \ra {\cal O}_{\bf P} \ra {\cal O}_X(n) \ra 0,$$
	which implies the sequence
	$$ 0 \ra \bigoplus\limits_{n \ge 0} \Gamma({\cal J}_X(n)) \ra
        \bigoplus\limits_{n \ge 0}
	\Gamma({\cal O}_{\bf P}(n)) \ra \bigoplus\limits_{n \ge 0} \Gamma({\bf P},
	{\cal O}_X(n)) \ra  \bigoplus\limits_{n \ge 0} H^1({\bf P}, {\cal J}_X(n)).$$
	Here the last term is trivial for sufficiently large $n$. The first
	three terms are equal correspondingly to $I, A$ and $S'$. It means that
	the homogenous components of $A/I$ and $S' \supset A/I$ are equal for
	sufficiently big degrees.

	By \cite{GD}[prop. 2.9.1]
	$$Proj(A/I) \cong Proj(S'),$$
	and combining everything together we get the statement of the lemma.
        The statement concerning the sheaf ${\cal F}$ can be proved along
        the same line.
\par\medskip
 Let us first  explain the Krichever correspondence for dimension 1.

 {\sc Definition 2}.
 $$
 \begin{array}{lll}
 {\cal M}_1 & := & \{ C, P, z, {\cal F}, e_P \}\\
 C & & \mbox{projective irreducible curve}~/k \\
 P \in C && \mbox{a smooth point} \\
 z && \mbox{formal local parameter at}~P \\
 {\cal F} && \mbox{torsion free rank}~$r$~\mbox{sheaf on}~C \\
 e_P && \mbox{a trivialization of}~{\cal F}~\mbox{at}~P
 \end{array}
 $$
 Independently, we have the field $K = k((z))$ of Laurent power series with
 filtration  $K(n) = z^nk[[z]]$. Let $K_1 : = K(0)$. If $V = k((z))^{\oplus r}$
 then $V(n) = K(n)^{\oplus r}$ and $V_1 : = V(0)$.

\begin{theorem}\cite{M}. There exists a canonical map
$$
\Phi_1 : {\cal M}_1 \ra \{\mbox{vector subspaces}~A \subset K, W \subset V \}
$$
such that
\begin{quotation}
i) the cohomology of complexes
 $$
A \oplus K_1 \ra K, ~W \oplus V_1 \ra V
 $$
are isomorphic to $H^{\cdot}(C, {\cal O}_C)$ and $H^{\cdot}(C, {\cal F})$,
respectively

ii) if $(A, W) \in \mbox{Im}~\Phi_1$ then $A \cdot A \subset A,
A \cdot W \subset W$,

iii) if $m,m^{\prime} \in {\cal M}_1$ and $\Phi_1(m) =
\Phi_1(m^{\prime})$ then $m$ is isomorphic to $m^{\prime}$
\end{quotation}
\end{theorem}

{\sc Proof}. If~$m = (C, P, z, {\cal F}, e_P) \in {\cal M}_1$ then we put
$$
A : = \Gamma(C - P, {\cal O}_C),
$$
$$
W: = \Gamma(C - P, {\cal F}).
$$
Also we have
$$
\cO_P = k[[z]],~K_P = k((z)),
$$
$$
{\cal F}_P  = {\cal O}_Pe_P = {\cal O}_P^{\oplus r},~\cF_P = \cO_P^{\oplus r}.
$$
This defines the point $\Phi_1(m) \in {\cal M}_1$. Indeed, for the subspace
$W$ we have the following canonical identifications
$$
\Gamma(C - P, {\cal F}) \subset {\cal F}_{\eta} \otimes_{{\cal O}_P} K_P
= \cF_P \otimes K_P = \cO_P^{\oplus r} \otimes K_P = k((z))^{\oplus r}.
$$
The same works for the subspace $A$.

The property ii) is
obvious, the property i) follows from the proposition 1. To get iii) let
us start with a point $\Phi_1(m) = (A, W)$. The standard valuation on $K$
gives us  increasing filtrations $A(n) = A \cap K(n)$ and
$W(n) = W \cap V(n)$ on the spaces $A$ and $W$. Then we have
$$
\begin{array}{ccl}
C - P & = & \mbox{Spec}(A),  \\
C  & = &  \mbox{Proj}(\oplus_n A(n)), \\
{\cal F} & = & \mbox{Proj}(\oplus_n W(n)),
\end{array}
$$
by lemma 6. Thus we can reconstruct the quintiple $m$ from the point
$\Phi_1(m)$.

{\sc Remark 4}. It is possible to replace the ground field $k$ in the
Krichever construction by an arbitrary scheme $S$, see \cite{Q}.

Now we move to the case of algebraic surfaces. The corresponding data
has the following

{\sc Definition 3}.
 $$
 \begin{array}{lll}
 {\cal M}_2 & := & \{ X, C, P, (z_1, z_2),  {\cal F}, e_P \}\\
 X & & \mbox{projective irreducible surface}~/k \\
 C \subset X & & \mbox{projective irreducible curve}~/k \\
 P \in C && \mbox{a smooth point on }~X~\mbox{and}~C \\
 z_1, z_2 && \mbox{formal local parameter at}~P~\mbox{such that} \\
        && (z_2 = 0) = C ~\mbox{near}~P \\
 {\cal F} && \mbox{torsion free rank}~$r$~\mbox{sheaf on}~X \\
 e_P && \mbox{a trivialization of}~{\cal F}~\mbox{at}~P
 \end{array}
 $$
Then we have
$$
\cO_{X,P}= k[[z_1, z_2]],~K_{P,C} = k((z_1))((z_2)),
$$
$$
\cF_P  = \cO_Pe_P = \cO_P^{\oplus r}.
$$
For the field $K = k((z_1))((z_2))$ we have the following filtrations and
subspaces:
$$
\begin{array}{lll}
K_{02} & = & k[[z_1]]((z_2)), \\
K_{12} & = & k((z_1))[[z_2]],   \\
K(n) & = & z_2^n K_{12}.
\end{array}
$$
Taking the direct sums we introduce the subspaces $V_{02}, V_{12},V(n)$ of
the space $V = K^{\oplus r}$.
\begin{theorem}. Let~$C$ be a hyperplane section on the surface $X$. Then
there exists a canonical map
$$
\Phi_2 : {\cal M}_2 \ra \{\mbox{vector subspaces}~B \subset K, W \subset V \}
$$
such that
\begin{quotation}

i) for all $n$ the complexes
 $$
 \frac{B \cap K(n)}{B \cap K(n+1)} \oplus \frac{K_{02} \cap K(n)}{K_{02} \cap
  K(n+1)} \ra \frac{K(n)}{K(n+1)}
$$
 $$
 \frac{W \cap V(n)}{W \cap V(n+1)} \oplus \frac{V_{02} \cap V(n)}{V_{02} \cap
  V(n+1)} \ra \frac{V(n)}{V(n+1)}
$$

are Fredholm of index $\chi(C, {\cal O}_C)  + nC.C$ and
$\chi(C, {\cal F}\vert_C)  + nC.C$, respectively

ii) the cohomology of  complexes
$$
(B \cap K_{02}) \oplus (B \cap K_{12}) \oplus (K_{02}) \cap K_{12})  \ra
B \oplus K_{02} \oplus K_{12} \ra K
$$
$$
(W \cap V_{02}) \oplus (W \cap V_{12}) \oplus (V_{02}) \cap V_{12})  \ra
W \oplus V_{02} \oplus V_{12} \ra V
$$
are isomorphic to $H^{\cdot}(X, {\cal O}_X)$ and $H^{\cdot}(X, {\cal F})$,
respectively

iii)if $(B, W) \in \mbox{Im}~\Phi_2$ then $B \cdot B \subset B,
B \cdot W \subset W$

iv) for all $n$ the map
$$
(C, P, z_1\vert_C, {\cal F}(nC)\vert_C, e_P(n)\vert_C) \mapsto
$$
$$
\mapsto (\frac{B \cap K(n)}{B \cap K(n+1)}  \subset  \frac{K(n)}{K(n+1)} =
k((z_1)),
$$
$$
\frac{W \cap V(n)}{W \cap V(n+1)} \subset  \frac{V(n)}{V(n+1)} =
k((z_1))^{\oplus r}~)
$$
coincides with the map $\Phi_1$.

v) let the sheaf ${\cal F}$ be locally free and the surface $X$ be Cohen-
Macaulay. If $m,m^{\prime}
\in {\cal M}_1$ and $\Phi_2(m) = \Phi_2(m^{\prime})$ then $m$ is isomorphic
to $m^{\prime}$
\end{quotation}
\end{theorem}

{\sc Proof}. If $ m = (X, C, P, (z_1, z_2),  {\cal F}, e_P) \in {\cal M}_2$ then
to define the map $\Phi_2$ we put
$$
\begin{array}{ccl}
B &  = & B_C({\cal O}_X), \\
W &  = & B_C({\cal F}), \\
\Phi_2(m) & = & (B, W).
\end{array}
$$
Since we have the local coordinates $z_{1,2}$ and the trivialization $e_P$
the subspaces $B$ and $W$ will belong to the space $k((z_1))((z_2))$ exactly
as in the case of dimension 1 considered above.

We note that our condition on the curve $C$ implies that $C$ ia  a Cartier
divisor and the surface $X - C$ is affine.

The property i) follows
from lemma 2, the property ii) follows from theorem 1. The property
iii) is trivial again, to get iv) one needs again to apply lemma 2 and
to get v) it is enough to use proposition 2 and lemma 7. They show that given
a point $(B, W) \in {\cal M}_2$ such that $(B, W) = \Phi_2(m)$ we can
reconstruct the data $m$ up to an isomorphism.

{\sc Remark 5}. The property v) of the theorem cannot be extended to the
arbitrary torsion free sheaves on $X$. We certainly cannot reconstruct
such sheaf if it is not locally free outside $C$. Indeed, if ${\cal F},
{\cal F}^{\prime}$ are two sheaves and there is a monomorphism
${\cal F}^{\prime} \ra {\cal F}$ such that ${\cal F}/{\cal F}^{\prime}$ has
support in $X - C$ then the restricted adelic complexes for the
sheaves ${\cal F}, {\cal F}^{\prime}$ are isomorphic.

{\sc Remark 6}. A definition of the map $\Phi_n$ for all $n$ was
suggested in \cite{O}. It has the properties that correspond to the
properties i) - v) of the theorem.

\end{document}